\documentclass[12pt]{amsart}
\usepackage{amssymb}
\usepackage{enumerate}

\usepackage{url}
\usepackage{hyperref}

\usepackage{vmargin}
\setpapersize{USletter}
\setmargrb{1in}{1in}{1in}{1in}

\newcommand{\Mustata}{Musta\c{t}\u{a}}

\newcommand{\C}{\ensuremath{\mathbb C}}
\newcommand{\R}{\ensuremath{\mathbb R}}

\newcommand{\J}{\ensuremath{\mathcal J}}
\newcommand{\A}{\ensuremath{\mathcal{A}}}
\newcommand{\G}{\ensuremath{\mathcal{G}}}
\newcommand{\Gmin}{\G_{\min}}
\renewcommand{\O}{\mathcal{O}}

\DeclareMathOperator{\codim}{codim}
\DeclareMathOperator{\mult}{mult}

\theoremstyle{definition}
\newtheorem{defn}{Definition}[section]
\newtheorem{example}[defn]{Example}
\newtheorem{remark}[defn]{Remark}

\theoremstyle{plain}
\newtheorem{lem}[defn]{Lemma}
\newtheorem{thm}[defn]{Theorem}

\title[A note on multiplier ideals of hyperplane arrangements]{A note on Musta\c{t}\u{a}'s computation of multiplier ideals of hyperplane arrangements}
\author{Zach Teitler}
\date{3/5/07}


\address{Department of Mathematics, SLU 10687, Hammond, LA 70402}
\email{zteitler@selu.edu}

\begin{document}

\bibliographystyle{plain}       

\begin{abstract}
In~\cite{mustata:hyperplane-arrangements}, M.~\Mustata\ used jet schemes
to compute the multiplier ideals
of reduced hyperplane arrangements.
We give a simpler proof using a log resolution
and generalize to non-reduced arrangements.
By applying the idea of wonderful models introduced by
De Concini--Procesi~\cite{MR1366622}, we also simplify the result.
Indeed, \Mustata's result expresses the multiplier ideal as an intersection,
and our result uses (generally) fewer terms in the intersection.
\end{abstract}

\maketitle

\section{Introduction}

For an ideal $I \subset \C[x_1,\dots,x_n]$,
regarded as an ideal on $\C^{n}$,
a \textbf{log resolution} of $I$ is a proper birational map
$f:X \to \C^n$, with $X$ smooth, such that the total transform
$I\cdot\O_X = \O_X(-F)$
is locally principal and
$F + \text{Exc}(f)$ is a divisor with normal crossings support.
Then for $\lambda \in \R$, $\lambda \geq 0$, the $\lambda$th
\textbf{multiplier ideal}
$\J(I^\lambda)$ is given by
\[
  \J(I^\lambda)
    = f_{*} \O_X( K_{X/\C^{n}} - \lfloor \lambda F \rfloor) .
\]
More details on multiplier ideals may be found in \cite{pag2}.

Let $\A$ be a hyperplane arrangement in $V \cong \C^n$.
For simplicity we assume $\A$ is central, that is, all hyperplanes
pass through the origin.
Suppose $\{H_1,\dots,H_r\}$ is the set of hyperplanes appearing in $\A$
and each $H_i$ is defined by the vanishing of a linear form $L_i$.
We allow the $H_i$ to have positive integer multiplicities:
for each $H_i$, let $m_i=\mult_{H_i}(\A)$ be the multiplicity of $H_i$ in $\A$.
The arrangement is reduced if every $m_i=1$.
Then the ideal $I=I(\A)=(L_1^{m_1}\cdots L_r^{m_r})$ defines $\A$.
The goal is to compute $\J(I^\lambda)$.

Let $L(\A)$ be the intersection lattice of $\A$, the set of all
intersections of hyperplanes in $\A$.
For $W \in L(\A)$, define the rank of $W$ to be $r(W) = \codim(W)$
and let
\[ s(W) = \mult_W(\A) = \sum_{W \subset H \in \A} \mult_H(\A) . \]
For a reduced arrangement, $s(W)$ is the number of hyperplanes of $\A$ containing $W$.
Let $L'(\A) = L(\A) \setminus \{ V \}$.
Then M.~\Mustata\ computes the multiplier ideals of a reduced hyperplane arrangement
in~\cite{mustata:hyperplane-arrangements}, obtaining the following result.
\begin{thm}\label{thm-mustata}
Let $\A$ be a reduced hyperplane arrangement, with defining ideal $I$.
Then for $\lambda \geq 0$,
\[
  \J( I^\lambda) = \bigcap_{W \in L'(\A)} I_W^{ \lfloor \lambda s(W) \rfloor - r(W) + 1}
\]
where $I_W$ is the ideal of $W$.
\end{thm}
This is proved using jet schemes.
Though it is not stated explicitly in~\cite{mustata:hyperplane-arrangements},
the method of jet schemes can treat the case of non-reduced arrangements.
(See also~\cite{saito:mis-of-locally-conical-divs} for a generalization
to a locally conical divisor along a stratification, proved by different methods.)
It is possible, however, to give a proof using simply a log resolution of the hyperplane arrangement,
as suggested by Remark 1.2 of~\cite{mustata:hyperplane-arrangements}.
Using the notion of building sets then allows us to simplify the result,
in the sense of replacing the intersection over $L'(\A)$ by an intersection
with possibly fewer terms.

We briefly recall the notion of building sets introduced by
De~Concini and Procesi~\cite[\textsection2.3]{MR1366622}.
We restrict to the special case of hyperplane arrangements.
(In~\cite{MR1366622}, arbitrary subspace arrangements are treated.
See also~\cite{MR2178326} for an expository account
and~\cite{yi-hu,math.AG/0611412} for generalizations.)
\begin{defn}
Let $\A$ be a hyperplane arrangement in $V$.
A \textbf{decomposition} of $C \in L'(\A)$
is a subset $\{U_1, \dots , U_k\} \subset L'(\A)$
such that $C = U_1 \cap \dots \cap U_k$,
transversally (that is, $\codim C = \codim U_1 + \dots + \codim U_k$);
and for any $C \subset B \in L'(\A)$,
we have each linear sum $B + U_i \in L(\A)$ and
$B = (B + U_1) \cap \dots \cap (B + U_k)$, again transversally.

A subset $\G \subset L'(\A)$ is a \textbf{building set} if
for every $C \in L'(\A)$, the minimal elements
$\{G_1, \dots, G_r \}$ of $\G$ containing $C$
give a decomposition of $C$.
\end{defn}

\begin{example}
\begin{enumerate}[(a)]
\item Each $C \in L'(\A)$ admits the trivial decomposition $\{C\}$.
Correspondingly, $L'(\A)$ is itself a building set.
\item 
An element in $L'(\A)$ is called \textbf{irreducible}
if it admits no non-trivial decomposition.
In particular, every hyperplane in $\A$ is irreducible.
De~Concini--Procesi show that
the set of irreducible elements forms a building set,
which we denote $\Gmin$.
It is containment-minimal in the sense that $\Gmin$ is contained
in every other building set~\cite[\textsection2.3]{MR1366622}.
\item
The \textbf{braid arrangement} $\mathcal{B}_n$ on $\C^n$
has hyperplanes $H_{ij}$ defined by $x_i=x_j$, for $1 \leq i < j \leq n$.
(Sometimes $\mathcal{B}_n$ is considered as an arrangement on $\C^{n-1}$
via quotienting out by the line $x_1=\dots=x_n$.)
The intersection lattice $L(\mathcal{B}_n)$ is isomorphic to
the lattice of partitions of $\{1,\dots,n\}$, ordered by reversed refinement.
For example, the subspace $W_{123 \mid 45} \in L(\mathcal{B}_n)$
is defined by the equations $x_1=x_2=x_3$ and $x_4=x_5$,
so it is the intersection $(H_{12}\cap H_{13}\cap H_{23}) \cap (H_{45})$.

Now, $W_{123} = H_{12} \cap H_{13}$ is a transversal intersection,
but $\{H_{12},H_{13}\}$ is not a decomposition of $W_{123}$.
Indeed, $W_{123} \subset H_{23}$, yet
\[ (H_{23}+H_{12}) \cap (H_{23}+H_{13}) = \C^3 \cap \C^3 \neq H_{23} \]

Let $p$ be a partition of $\{1,\dots,n\}$ with blocks $b_1, \dots, b_k$ of size greater than $1$.
Then $W_p$ admits the decomposition $\{ W_{b_1}, \dots, W_{b_k} \}$.

Conversely, if $p=b_1$ has only one block of size greater than $1$,
then $W_p$ is irreducible.
Such partitions are called \textbf{modular}.
It follows that in the braid arrangement $\mathcal{B}_n$,
the minimal building set $\Gmin$ consists of $W_p$ with $p$ modular.
For $n \gg 0$,
\[ \# \Gmin = 2^n - n - 1 \ll \#L(\mathcal{B}_n) \]
In fact, the numbers $\# L(\mathcal{B}_n)$, called Bell numbers~\cite{bellnumbers},
are super-exponential.
For example, with $n=10$,
$\#\Gmin=1{,}013$,
$\#L(\mathcal{B}_{10})=115{,}975$;
with $n=20$, $\#\Gmin/\#L(\mathcal{B}_{20}) \approx 2.03 \cdot 10^{-8}$.
\end{enumerate}
\end{example}

\begin{thm}\label{thm:buildingsets}
Let $\A$ be a hyperplane arrangement (not necessarily reduced) with ideal $I$.
Let $\G \subset L'(\A)$ be a building set.
Then for $\lambda \geq 0$,
\begin{equation}\label{eqn:buildingsets}
  \J( I^\lambda) = \bigcap_{W \in \G} I_W^{ \lfloor \lambda s(W) \rfloor - r(W) + 1}
\end{equation}
\end{thm}
Theorem~\ref{thm-mustata} is the case $\G=L'(\A)$,
and $\A$ reduced.
The minimal building set $\Gmin$ gives the version of~(\ref{eqn:buildingsets})
with the fewest terms in the intersection.
The example of the braid arrangement shows this can have dramatically fewer terms.

\section{Log resolution}

Let $\A$ be a hyperplane arrangement in $V$, not necessarily reduced,
and let $\G \subset L'(\A)$ be a building set.
Let $\G = G_0 \cup G_1 \cup \dots \cup G_n$,
where $G_i = \{ \, W \in \G \mid \dim(W) = i \, \}$.
We blow up the space $V$ iteratively: First blow up $G_0$;
then blow up the proper transforms of all subspaces in $G_1$; and so on.
At each stage, the spaces to be blown up are disjoint because their intersections
have been blown up already at an earlier stage.

We denote this space by $V_{\G}$, with blowdown $f: V_{\G} \to V$.
It is shown in~\cite{MR1366622} that the set-theoretic preimage
$f^{-1}(\A)$ is a divisor with simple normal crossings support.
The following lemma will show $f$ is a log resolution of $I=I(\A)$.

Each $W \in \G$ is dominated by
a unique prime divisor $E_W$ in $V_{\G}$.
For $W$ not a hyperplane in $V$, $E_W$ is $f$-exceptional.
It is the proper transform in $V_{\G}$
of the exceptional divisor produced by blowing up
(the proper transform of) $W$ in an earlier stage.
For $W=H_i$ a hyperplane in $\A$, blowing up
(the proper transform of) $W$ is the identity map.
In this case $E_W$ is just the proper transform of $W$,
so it is not $f$-exceptional.

\begin{lem}\label{lem:resolution}
Let $\A \subset V$ be a possibly non-reduced hyperplane arrangement
with ideal $I$.
Let $\G \subset L(\A)$ be a building set.
The map $f: V_{\G} \to V$
is a log resolution of $I$.
For $W \in \G$, let $E_W \subset V_{\G}$ be the prime
divisor dominating $W$.
The relative canonical divisor is
\[ K_{V_{\G}/V} = \sum_{W \in \G} (r(W)-1)E_W, \]
where $r(W) = \codim_V(W)$.
The pullback $f^{*}(\A)$ is
\[ f^{*}(\A) = \sum_{W \in \G} s(W) E_W, \]
where as above $s(W) = \mult_W(\A) = \sum_{W \subset H \in \A} \mult_H(\A)$.
\end{lem}

\begin{proof}
The pullback $f^{*}I$ is (the ideal of) a divisor supported along the set-theoretic preimage
$f^{-1}(\A)$, which is a divisor with normal crossings support.
The exceptional locus of $f$ also has support contained in $f^{-1}(\A)$.
This shows $f$ is a log resolution of $I$.

For the description of the relative canonical divisor, see~\cite[Exer.~II.8.5(b)]{hartshorne}.

For $H_i \in \A$, let $H_i$ be defined by the linear form $L_i$.
Then $f^{*}(L_i)$ vanishes along $E_W$ to order $0$ or $1$,
according as $W \not\subset H_i$ or $W \subset H_i$.
It follows that $f^{*}I = f^{*}(L_1^{m_1}\cdots L_r^{m_r})$ vanishes along $E_W$
to order $s(W)$, as claimed.
\end{proof}

\begin{remark}
More generally,~\cite{MR1366622} considers (linear) subspace arrangements.
In this more general setting,
the same idea of iteratively blowing up along a building set will give what 
the authors call a ``wonderful model'' of the subspace arrangement,
meaning a proper, birational map which is an isomorphism over the complement of the support
of the arrangement and such that the set-theoretic preimage of the arrangement
is a divisor with normal crossings support.
This is not always a log resolution, however, since outside the case of hyperplane arrangements
there may arise embedded components in the pullback of the ideal of the arrangement~\cite{zct:mla}.
\end{remark}

\section{Multiplier ideals}

We prove Theorem~\ref{thm:buildingsets}.
The key is the following lemma.

\begin{lem}\label{lem:pushforward}
With notation as in Lemma~\ref{lem:resolution}, for $p \geq 0$,
\[ f_{*} \O_{V_{\G}}(-p E_W) = I_W^p . \]
\end{lem}
\begin{proof}
Let $f: V_{\G} \to V$ be decomposed into stages of blowing up:
\[ V_{\G} = V_{n} \to V_{n-1} \to \dots \to V_0 = V , \]
where $V_{i+1} \to V_{i}$ is the blowing-up of the (proper transforms of the) subspaces
in $\G$ of dimension $i$.
Let $d=\dim(W)$ and consider 
\[ V_{\G} \overset{c}{\to} V_{d+1} \overset{b}{\to} V_d \overset{a}{\to} V . \]
We denote the proper transform of $W$ in $V_d$ by $W'$.
We denote the irreducible exceptional divisor in $V_{d+1}$ over $W'$ simply by $E$.
Then $E_W$ is the total transform of $E$.
It follows $c_{*} \O_{V_{\G}}(-p E_W) = \O_{V_{d+1}}(-p E)$.
Since $W'$ is smooth, it follows immediately that
$b_{*} \O_{V_d}(-p E) = I_{W'}^p$,
where $I_{W'}$ is the ideal sheaf of $W'$.
And since $W'$ is the proper transform of $W$,
$a_{*} I_{W'}^p = I_W^p$.
\end{proof}

\begin{proof}[Proof of Theorem~\ref{thm:buildingsets}]
Taking the log resolution of Lemma~\ref{lem:resolution}, we have
\[
  \begin{split}
  \J(I^\lambda) &= f_* \O_{V_{\G}} ( K_{V_{\G}/V} - \lfloor \lambda F \rfloor ) \\
    &= f_* \O_{V_{\G}} \Big( \sum_{W \in \G} ( r(W)-1-\lfloor \lambda s(W) \rfloor ) E_W  \Big) \\
    &= \bigcap_{W \in \G} f_* \O_{V_{\G}} \big( ( r(W) - 1 - \lfloor \lambda s(W) \rfloor ) E_W \big) \\
    &= \bigcap_{W \in \G} I_W^{ \lfloor \lambda s(W) \rfloor - r(W) + 1 }
  \end{split}
\]
with the last equality following from Lemma~\ref{lem:pushforward}.
\end{proof}

\begin{remark}
We can slightly refine two corollaries of~\cite{mustata:hyperplane-arrangements}.
We have from Corollary~0.2 that (using notation from above) the support of $\J(I^{\lambda})$
is the union of those $W \in \Gmin$ with $\lambda \geq r(W)/s(W)$.
From Corollary~0.3 we see that the log canonical threshold of $I$
is
\[ \text{lct}(I) = \min_{W \in \Gmin} \frac{ s(W) }{ r(W) } . \]
In each case we have replaced the condition $W \in L'(\A)$ with $W \in \Gmin$,
and removed the condition $\A$ be reduced.

Example 2.3 of~\cite{mustata:hyperplane-arrangements} (concerning set-theoretic jumping
numbers) admits a similar refinement.
\end{remark}



\section*{Acknowledgments}
The author thanks Mircea~\Mustata, Hal~Schenck,
and the referee for a number of very helpful comments.

\bibliography{../biblio}   

\end{document}